\documentclass[11pt]{article}
\def\version{29.3.'19}\def\users{}  %

\def\users{final-layout}  

\usepackage{amsmath}
\usepackage{paralist}
\usepackage{cite}  
\usepackage{graphics} 
  \usepackage{epsfig} 
\usepackage[colorlinks=true]{hyperref}
\hypersetup{urlcolor=blue, citecolor=red}

\newtheorem{proposition}{Proposition}

\newtheorem{definition}{Definition}
\newtheorem{remark}{Remark}

\usepackage{amssymb}
\usepackage{pst-all}

\usepackage{psfrag}
\newcounter{myfigure}
\newenvironment{my-picture}[3]{\refstepcounter{myfigure}\label{#3}\setlength{\unitlength}{\textwidth}\begin{picture}(#1,#2)}{\end{picture}}
\usepackage{ifthen}
\ifthenelse{\equal{\users}{final-layout}}{}
{
\usepackage{fancyhdr}
\pagestyle{fancy}
\headheight=28pt
\rhead{\color{green}
Cahn-Hilliard actual...\\
by T.Roub\'\i\v{c}ek}
\chead{}
\lhead{\color{green}Version \version, file:\jobname.tex\\compiled:
\number\day.\number\month.\number\year\ at
\the\hour:\ifnum\minute<10 0\fi\the\minute\ h }

\usepackage[notcite,notref,color]{showkeys}
\definecolor{labelkey}{rgb}{1.,.2,0.}

\newcount\hour \newcount\minute
\hour=\time
\divide \hour by 60
\minute=\time
\loop \ifnum \minute > 59 \advance \minute by -60 \repeat

}

\usepackage[normalem]{ulem}
\definecolor{brown}{rgb}{0.5,0,0}
\ifthenelse{\equal{\users}{final-layout}}{

	\newcommand{\COMMENT}[1]{}
	\newcommand{\DELETE}[1]{}

        \newcommand{\REM}[1]{\marginpar{\bfseries\tiny{\color{blue}}}}
         }
{
 
 \newcommand{\COMMENT}[1]{{\color{red}\uuline{#1}\color{black}}}
 \newcommand{\DELETE}[1]{{\color{brown}\sout{#1}\color{black}}}

 \newcommand{\REM}[1]{\marginpar{\bfseries\tiny{\color{blue}#1}}}
}

\usepackage{mathrsfs}   

\newcommand{\PP}{\varPi}
\newcommand{\eq}[1]{(\ref{#1})}
\newcommand{\calE}{\mathcal E}
\newcommand{\bbH}{\mathbb H}

\newcommand{\bbM}{\mathbb M}

\newcommand{\Vdots}{\mathchoice{\,\vdots\,}{\:\begin{minipage}[c]{.1em}\vspace*{-.4em}$^{\vdots}$\end{minipage}\;}
{\:\tiny\vdots\:}{\:\tiny\vdots\:}}
\newcommand{\In}{\!\in\!}
\newcommand{\pl}{\partial}
\newcommand{\plF}{\partial_{\!F}^{}}
\renewcommand{\d}{\mathrm d}  
\newcommand{\nablaS}{\nabla_{\scriptscriptstyle\textrm{\hspace*{-.3em}S}}^{}}
\newcommand{\nablaSS}{\nabla_{\scriptscriptstyle\textrm{\hspace*{-.3em}S}}^2}
\newcommand{\divS}{\mathrm{div}_{\scriptscriptstyle\textrm{\hspace*{-.1em}S}}^{}}
\newcommand{\divSS}{\mathrm{div}_{\scriptscriptstyle\textrm{\hspace*{-.1em}S}}^2}

\newcommand{\Colon}{\!\colon\!}
\newcommand{\Cdot}{\!\cdot\!}
\newcommand\DT[1]{\mathchoice
                 {{\buildrel{\hspace*{.1em}\text{\LARGE.}}\over{#1}}}
                 {{\buildrel{\hspace*{.1em}\text{\Large.}}\over{#1}}}
                 {{\buildrel{\hspace*{.1em}\text{\large.}}\over{#1}}}
                 {{\buildrel{\hspace*{.1em}\text{\large.}}\over{#1}}}}
\newcommand\DDT[1]{\mathchoice
   {{\buildrel{\hspace*{.13em}\text{\LARGE.\hspace*{-.13em}.}}\over{#1}}}
   {{\buildrel{\hspace*{.1em}\text{\Large.\hspace*{-.1em}.}}\over{#1}}}
   {{\buildrel{\hspace*{.1em}\text{\large.\hspace*{-.1em}.}}\over{#1}}}
   {{\buildrel{\hspace*{.1em}\text{\large.\hspace*{-.1em}.}}\over{#1}}}}

\begin{document}
\begin{sloppypar}

\begin{center}
{\LARGE\bf Cahn-Hilliard equation with capillarity\\[.3em]in actual deforming configurations} \footnote{
This research has been partially supported from the grants
17-04301S
(especially as far as the focus on the dissipative evolution of internal
variables) and 
19-04956S
(especially as far as the focus on the dynamic and nonlinear behaviour)
of Czech Science Foundation,
and the FWF/CSF project 19-29646L
(especially as far as the focuse on the large strains in materials science),
 and also from the 
institutional support RVO:\,61388998 (\v CR).}
\end{center}

\def\SS{S}
\def\SSel{S_{\!\text{\sc el}}^{}}
\def\SSin{S_{\!\text{\sc in}}^{}}
\def\SSvi{S_{\!\text{\sc vi}}^{}}
\def\GG{G}
\def\R{\mathbb R}
\def\N{\mathbb N}
\def\GDir{\varGamma_{\rm Dir}^{}}
\def\Rsym{\R_{\rm sym}^{d\times d}}


\centerline{\scshape T.~Roub\'\i\v cek}


\medskip



\begin{center}
{\footnotesize
Institute of Thermomechanics, \,Czech Acad.\ Sci.,\\[-.2em]
Dolej\v skova 5, CZ-18200 Praha 8, Czech Republic 
}
\end{center}

\bigskip\bigskip

\centerline{\Large\it Dedicated to Alexander Mielke on the occasion of his sixtieth birthday.}

\bigskip\bigskip

\begin{quote}{\normalfont\fontsize{8}{10}\selectfont
    {\bfseries Abstract.}
    The diffusion driven by the gradient of the chemical potential
    (by the Fick/Darcy law) in deforming continua at large strains is
    formulated in the reference configuration with both the Fick/Darcy law
    and the capillarity gradient term considered at the actual configurations
    deforming in time. Static situations are analysed by the direct method.
    Evolution (dynamical) problems are treated by the Galerkin method, the
    actual capillarity giving rise to various new terms as e.g.\ the
    Korteweg-like stress and analytical difficulties related to them.
    Some other models (namely plasticity at small elastic strains or
    damage) with gradients at actual configuration allow for similar models
    and analysis.
\par

{AMS Subj. Classification:
 35Q74, 
37N15, 
65M60, 
74A30, 
74G25, 
74H20, 
76S05. 
}


  }
\end{quote}

\bigskip

\section{Introduction}\label{sect-intro}


\def\bmy{y}
\def\bmF{F}
\def\wt{\widetilde}

This paper addresses the Cahn-Hilliard model \cite{CahHil58FEUS} for diffusion
with capillarity (i.e.\ concentration gradient involved in the stored
energy) in deformable media, which is sometimes also called the Cahn-Larch\'{e}
model \cite{LarCah82ESSD}.
It is usually considered in mathematical literature
mostly at small strains, cf.\ e.g.\ \cite{BCDGSS02MPSBA,Garc03CHSE,PawZaj08WSCH,RouTom14THSM} and
\cite[Ch.4,5,7]{HeiKra14PSCD} and references therein.
If at large strains, then it is usually considered
with the concentration gradient in the reference (undeformed)
configuration, 
as used also in most of engineering references,
cf.\ e.g.\ \cite{ArSaRa16SACC,DLReAn14CHTP}, and as is relevant in some
applications.

Yet, in some other applications, the gradients in the actual space
(deforming) configurations seem more natural.
As mentioned in \cite{HGOD17FPLS} in the context of Allen-Cahn equation
about these gradients,
``their spatial counterparts could have also been used, this would
lead to cumbersome contributions (via pull-back and push-forward operations)''.
Nevertheless, when the transport by Fick/Darcy law is considered in
actual space deforming configuration as e.g.\ in
\cite{Anan12CHTT,KruRou19MMCM,RouSte18TEPR,RouTom13TMPE},
it rather mis-conceptual to involve gradient of concentration in the material
configuration. And indeed, sometimes the Cahn-Hilliard 
with the concentration gradient in the actual configuration 
can be found in engineering literature \cite{DalMie15CECM,HonWan13PFMS,MiMaUl14FNEM}
but, of course, without any analysis. A fairly general model has been
scrutinised in \cite{Pawl06TCCH} but without caring about non-selfinterpenetration
(and thus not considering possible singularities of the stored energy) and 
also without analysis as far as the existence of weak solutions concerns.

We will confine ourselves to a single-component flow. A generalization for
a multicomponent flow possibly also with mutual reactions between particular
components is interesting and seems possible, in particular when the 
gradient structure is kept, cf.\ \cite{Miel11GSRD,Miel133TMER}.

The main goal of this article is to perform a rigorous analysis as far as
the existence of weak solutions for the diffusion in poro-elastic-dynamic
model both with Fick/Darcy law and the capillarity gradient term in the
actual deforming configuration. 

First, in Section~\ref{sect-static}, we will specify the stored energy and
present the static problem, exploiting the 2nd-grade  nonsimple material
concept, and perform the analysis as far as the existence of the
solution based on minimum-energy principle. Then, in Section~\ref{sect-dynam},
we will slightly modify the problem (by simplifying the constraints and
using 3rd-grade nonsimple material concept) and formulate the evolutionary
variant, involving also inertial effects, which allows for modelling
elastic waves interacting with diffusion equation, e.g.\ waves whose attenuation
and dispersion can be influenced by the content of diffusant. For the
analysis, we use the Galerkin approximation which can keep
the approximate solutions well away from the singularity of the
stored energy at deformation gradients with non-positive determinants.
Eventually, in Section~\ref{sect-remarks}, we present various
other application of the presented mathematical techniques to gradient
theories for some other internal variables.

\section{Static problem}\label{sect-static}

Before treating the evolution problems, let us begin with static situations.
The equilibria of poroelastic or swelling-exhibiting materials are assumed to
be governed by the energy minimization and lead to interesting mathematical
problems. 

As usual in continuum mechanics of solids, we consider the Lagrangian
formulation with $\varOmega\subset\R^d$ a fixed reference domain.
The state variables are the deformation $\bmy:\varOmega\to\R^d$ and the
concentration $\zeta:\varOmega\to\R$. 
The basic ingredient for the model is the stored energy, here considered as
\begin{align}
  \calE(y,\zeta):=&\int_\varOmega\phi(\nabla\bmy,\zeta)
  +\frac{\kappa}2\big|(\nabla\bmy)^{-\top}\nabla\zeta\big|^2
  +\frac1p|\nabla^2\bmy|^{p-2}\nabla^2\bmy\Vdots{\mathbb H}\Vdots\nabla^2\bmy
  \,\d x
\label{stored-static}\end{align}
with $\kappa>0$ a capillarity coefficient and with
$\bbH$ 
a
(pressumably small) regularizing 6th-order symmetric positive definite
tensor, $p>d$.
The so-called  2nd-grade nonsimple-material (or couple-stress) concept
\cite{FriGur06TBBC,Toup62EMCS} has been applied, leading 
to the bending-like energy contribution due to the $\mathbb H$-term,
involving second-order 
deformation gradients (=\,first-order strain gradients).
In dynamical situations, this may offer a suitable 
tool to model a dispersion. Beside such mechanical motivation, the main
mathematical advantage of the nonsimple-material concept is that higher-order
deformation gradients bring additional regularity of deformations and
also compactness of the set of 
admissible deformations in a stronger topology. Moreover, there the stored 
energy can be even convex in the highest derivatives of the deformation, which 
is helpful in proving existence of minimizers.

Let us emphasize that the capillarity in \eq{stored-static} is considered in the
actual configuration, being pulled back into the reference configuration by a
vectorial pushforward $(\nabla\bmy)^{-\top}$.
Thus the determinant of $\nabla\bmy$ is to be kept away from zero
to have $(\nabla\bmy)^{-\top}$ under control, which needs involvement 
of the  $\mathbb{H}$-term. This also allows for most general,
not necessarily polyconvex stored energies.

In the static situations we are addressing in this section, all dissipation
processes vanishes, i.e.\ here in particular all diffusive processes vanish. 
Here it means that the gradient of the chemical potential vanishes
on $\varOmega$. When assuming $\varOmega$ connected, this further 
leads to that $\mu$ is constant, cf.\ also Remark~\ref{rem-mu-const}.
Let us denote this constant by $\bar\mu$. 

A variationally interesting situation is that the poroelastic body is
completely isolated on its boundary. 
%
It is then natural to prescribe the total 
amount of diffusant
\begin{align}\label{swelling-small-constraint}
\int_\varOmega \zeta\,\d x=Z^{}_{\rm total}\qquad\text{with\ \ \ }Z^{}_{\rm total}\ge0
\text{\ \ \ given}.
\end{align}


We will use the standard notation concerning the Lebesgue and the Sobolev
spaces, namely $L^p(\varOmega;\R^n)$ for Lebesgue measurable functions $\varOmega\to\R^n$ whose Euclidean norm is integrable with $p$-power, and
$W^{k,p}(\varOmega;\R^n)$ for functions from $L^p(\varOmega;\R^n)$ whose
all derivative up to the order $k$ have their Euclidean norm integrable with
$p$-power. We also write briefly $H^k=W^{k,2}$. Moreover,
${\rm GL}^+(d)$ denotes the general linear group of
orientation-preserving mappings $\R^d\to\R^d$, i.e.\ the subset of
$\R^{d\times d}$ of nonsingular matrices with a positive determinant,
while ${\rm SO}(d)$  denotes the special orthogonal group, i.e.\
the set $\{A\In\R^{d\times d};\ A^\top\!A=AA^\top=\mathbb I,\ \det A=1\}$.
 
We require that admissible deformations of the material are orientation
preserving and injective almost everywhere in $\varOmega$. The attribute will
be ensured by the Ciarlet-Ne\v{c}as condition \cite{CiaNec87ISCN}. We also
assume that the elastic body is fixed on a part of its boundary by a Dirichlet
condition. Altogether, we are left with the following problem: 
\begin{align}\label{poroelastic-largestrain}
\left.\begin{array}{ll}
\text{Minimize }&\displaystyle{J(y,\zeta):={\mathcal E}(y,\zeta)
  -\int_\varOmega
  f\Cdot y   \,\d x}\ \ 
\\[-.2em]\text{subject to }&
\displaystyle{\int_\varOmega\det\nabla y
\,\d x\le {\rm meas}_d(y(\varOmega))
\ \ \text{ and }\ \int_\varOmega \zeta\,\d x=Z_{\rm total}\,,\ \ \ }
\\[.8em]&
\det\nabla y>0\ \text{ and }\ c\ge 0\ \text{\ \ a.e.\ on }\varOmega,
\\[.3em]&y|_{\GDir}^{}=y_{\rm D}^{},\ \ \ y\In W^{2,p}(\varOmega;\R^d)\
\text{ and }\ \ \zeta\In H^1(\varOmega).
\end{array}\right\}
\end{align}
The physically motivated assumptions on the stored-energy density are
\begin{subequations}\label{ass-growrt-phi-dam}\begin{align}
    &\varphi:{\rm GL}^+(d)\times\R\to\R^+\ \
    \text{ continuously differentiable and
}
    \\&\label{frame-indif}\forall R\In{\rm SO}(d),\ \ z\In\R:\ \
    \varphi(F,z)=\varphi(RF,z)\,,
    \\&\nonumber
\exists\epsilon>0\ \ \forall F\In{\rm GL}^+(d),\ c\In\R:\ \ 
\\[-.3em]&\qquad\varphi(F,z)\begin{cases}
  \ge\displaystyle{\frac\epsilon{(\det\,F)^q}+\epsilon|z|^2}\!\!&\text{for }\det\,F>0
  \ \text{ with }\ \displaystyle{q>\frac{pd}{p{-}d}},\ \ p>d\,,\\[-.2em]
  =+\infty&\text{for }\det\,F\ge0\,.\end{cases} 
\label{growth}\end{align}\end{subequations}
   The assumption \eq{frame-indif} is the frame-indifference, while 
   \eq{growth} grants local non-selfinterpenetration and even allows to keep the
   deformation gradient ``uniformly'' invertible due to \cite{HeaKro09IWSS}.
    
\begin{proposition}
  Let $\varphi$ satisfy \eq{ass-growrt-phi-dam}, $\mathbb H$ be symmetric
  positive definite, 
  $y_{\rm D}^{}\In W^{2-1/p,p}(\GDir;\R^d)$, $f\In L^1(\varOmega)$, and
  \eq{poroelastic-largestrain} be feasible in the sense that its
  constraints are satisfied for at least one $(y_0,\zeta_0)$ with
$\inf_\varOmega\det(\nabla y_0)>0$.
  Then \eq{poroelastic-largestrain} has a solution $(y,\zeta)\in
  W^{2,p}(\varOmega;\R^d)\times H^1(\varOmega)$ such that
  $\inf_\varOmega\det(\nabla y)>0$.
\end{proposition}
  
\noindent{\it Proof.}  The assumed feasibility ensures existence of some $(y_0,\zeta_0)$ which is
  compatible with the constraints and which makes
  the functional \eq{stored-static} finite. By the Healey-Kr\"omer theorem 
  \cite{HeaKro09IWSS} and by the assumption \eq{growth}, there exists
  $\varepsilon>0$ such that $\det\nabla y\ge\varepsilon$ for any
  $(y,\zeta)$ from the respective level-set of $J$ for which
  $J(y,\zeta)\le J(y_0,\zeta_0)$.

  This makes the functional $J$ weakly lower-semicontinuous on this level-set.
  For any infimizing sequence $\{(y_k,\zeta_k)\}_{k\in\N}$, one can take a
  subsequence weakly converging in $W^{2,p}(\varOmega;\R^d)\times H^1(\varOmega)$.
  By compact embedding, $(\nabla y_k,\zeta_k)$ converges strongly in
  $L^\infty(\varOmega;\R^d)\times L^2(\varOmega)$ which makes the functional
  $(y,\zeta)\mapsto\int_\varOmega\varphi(\nabla y,\zeta)\,\d x$ lower-semicontinuous
  by the Fatou lemma.
  Moreover, also
  $$
  (\nabla y_k)^{-\top}=\frac{{\rm Cof}\nabla y_k}{\det\nabla y_k}
    \to\frac{{\rm Cof}\nabla y}{\det\nabla y}
    =(\nabla y)^{-\top}\quad\text{ strongly in }\ L^\infty(\varOmega;\R^{d\times d})\,,
      $$
      so that the functional $(y,\zeta)\mapsto
      \int_\varOmega\frac{\kappa}2\big|(\nabla\bmy)^{-\top}\nabla\zeta\big|^2\,\d x$
      is weakly lower-semi\-continuous.

      The existence of a minimizer follows then by the direct-method arguments.
$\hfill\Box$

      \medskip
      
The gradient terms can be omitted when $\varphi$ is so-called cross-polyconvex,
cf.\ \cite[Sect.\,3.6.1]{KruRou19MMCM}. 
In some cases, at least the $\mathscr{H}$-term can be omitted even if $\varphi$ is not cross-polyconvex but then the capillarity term is to be considered 
in the reference configuration, cf.\ \cite[Sect.\,3.6.2]{KruRou19MMCM}.

\begin{remark}[\emph{Constancy of chemical potential.}]\label{rem-mu-const}\upshape
From the optimality conditions for a solution $(y,\zeta)$ to
\eq{poroelastic-largestrain}, in particular involving the partial differential
with respect to $\zeta$, one can read formally (if $\varphi$ is
suitably smooth and the constraint $\zeta\ge0$ is not active) that there is
a Lagrange multiplier $\bar\mu\in\R$ to the constraint
$\int_\varOmega\zeta\,\d x-Z_{\rm total}=0$ and
$\pl_\zeta(\calE(y,\zeta)-\mu(\int_\varOmega\zeta\,\d x-Z_{\rm total}))=0$, i.e.
$$
\bar\mu
=\pl_z\varphi(\nabla y,\zeta)
-{\rm div}\big(\kappa(\nabla y)^{-1}(\nabla y)^{-\top}\nabla\zeta\big)
$$
on $\varOmega$. This multiplier is in the position of chemical potential.
\end{remark}

\begin{remark}[\emph{Steady-state problems.}]\upshape
  An interesting generalization would be towards steady-state problems where
  diffusion flux (being constant in time) does not necessarily vanish, cf.
\eq{CH-evol} below with all time-derivative omitted.
  Existence for such generalization seems
  open, however. Some results are available only at small strains by using the
  Schauder fixed-point theorem, cf.\ \cite{Roub17VMSS}. On the other
  hand, sometimes some self-induced oscillations in such porous media
  (polymer gels) are 
  observed, cf.\ e.g.\ \cite{YaBa07TCMS,YSYB12CDBH}, which may indicate that
  there might be even some physical reasons for nonexistence of steady-state
  solutions.
\end{remark}


\section{Dynamical problem}\label{sect-dynam}

Our main goal is to formulate evolution governed by the stored energy
from Sect.\,\ref{sect-static} and execute its analysis. We will focus
to dynamical problems, i.e.\ involving inertia. In contrast to
static situations, variational formulations (based now, instead of
minimal-energy principle, on the Hamiltonian
variational principle extended for nonconservative systems) do not
seem fitted with applications of direct methods. Instead, we will
use formulation in terms of conventional partial-differential equations with
corresponding boundary conditions.

For this reason, we need to adopt two compromising modifications
of the static problem. First, we will ignore the Ciarlet-Ne\v cas global
non-selfinterpenetration condition while keeping only the local
non-selfinterpenetration $\det(\nabla y)>0$ which is anyhow needed to
keep control under the pulled-back concentration gradient and the pull-backed
mobility gradient.
Second, we will also ignore the constraint $\zeta\ge0$; in fact, this is
an often accepted modelling simplification relying that the mobility
of the diffusant is very small and the stored energy very large if
concentration approaches zero. Third, 
we need to have the regularizing $\mathscr{H}$-term quadratic
so that the resulted nonlinear hyperbolic problem is linear in the
highest-order terms, which needs to involve (possibly fractional)
derivative of the deformation gradient of the order higher than
$1+d/2$. This is inevitably rather technical; for the fractional-gradient
and thus the concept of a nonlocal nonsimple material see e.g.\
\cite{KruRou19MMCM}. Here we take the option of the 3rd-grade nonsimple
material like in \cite{AgiLaz09CBLL,Mind65SGSS}, considering the stored energy
\begin{subequations}\label{functionals}\begin{align}\label{stored-dynam}
\calE(y,\zeta)=\int_\varOmega\varphi(\nabla\bmy,\zeta)
  +\frac{\kappa}2|(\nabla\bmy)^{-\top}\nabla\zeta|^2
+\frac12\nabla^3\bmy\Vdots\bbH\Vdots\nabla^3\bmy \,\d x 
\end{align}
with $\bbH$ some symmetric positive definite 8th-order tensor. 
%

The other ingredients in building the evolution model are the kinetic
energy
\begin{align}\label{T}
\mathcal T(\DT y):=\int_\varOmega\frac\varrho2|\DT y|^2\,\d x
\end{align}
with $\varrho>0$ the mass density and the dot denoting the time derivative,
and the (Rayleigh's pseudo)potential of dissipative forces related with
diffusion:
\begin{align}\label{diss-pot}
  \mathcal R(y,\zeta;\DT\zeta):=
  \int_\varOmega\frac12|\mathfrak{M}^{1/2}(\nabla y,\zeta)
  \nabla\Delta_{\mathfrak{M}(\nabla y,\zeta)}^{-1}\DT\zeta|^2\,\d x
+\int_\varGamma\frac12\alpha(\Delta_{\mathfrak{M}(\nabla y,\zeta)}^{-1}\DT\zeta)^2\,\d S,
\end{align}
where $\alpha>0$ is a phenomenological permeability coefficient of the boundary
and where $\Delta_{\mathfrak{M}}^{-1}:(r,\mu_{\rm ext}^{})\mapsto\mu$ is the linear
operator $H^1(\varOmega)^*\to H^1(\varOmega)$ defined by $\mu:=\,$the weak solution
to the equation ${\rm div}(\mathfrak{M}\nabla\mu)=r$ with the
Robin boundary conditions $\mathfrak{M}\nabla\mu\Cdot\vec{n}+\alpha\mu=\alpha\mu_{\rm ext}^{}$;
in the case $\mu_{\rm ext}^{}=0$ cf.\ \cite[Sect.\,5.2.6]{MieRou15RIST}. Eventually, we
consider the mechanical load $F$ determined by the bulk force $f$ and the surface
load $g$ as by the external chemical potential $\mu_{\rm ext}$ by
\begin{align}\label{load}
F(t,\wt y,\wt\zeta)\cong\big\langle F(t),(\wt y,\wt\zeta)\big\rangle=
\int_\varOmega f(t)\Cdot\wt y\,\d x+\int_\Sigma g(t)\Cdot\wt y+\alpha\mu_{\rm ext}^{}\wt\zeta\,\d S\,.
\end{align}
\end{subequations}
Let us note that the dissipation potential \eq{diss-pot} is
nonlocal. A natural requirement for thermodynamical consistency
(i.e.\ non-negative entropy production) is that
$\mathfrak{M}$ is positive semidefinite, so that its square root
$\mathfrak{M}^{1/2}$ occurring in \eq{diss-pot} has a good sense.

The notation
$\mathfrak{M}:\varOmega\times\R^{d\times d}\times\R\to\R_{\rm sym}^{d\times d}$
stands for the mobility tensor which occurs in the
generalized Fick law making the flux of the diffusant proportional
to the gradient of the chemical potential denoted by $\mu$.  
Consistently with the capillarity
in the actual configuration pulled-back, a reasonable
modeling concept that this Fick law
(in particular covering also Darcy law) is considered in the actual deforming
(time-dependent) configuration, and is then to be pulled back into the
fixed reference configuration. 
The {\it transformed} {\it Fick law} 
(i.e.\ {\it pulled back}) uses the matrix of mobility coefficients as 
\begin{align}\label{M-pull-back}
  \mathfrak{M}(x,F,z)
  :=
\frac{({\rm Cof }F^\top)\bbM(x,z){\rm Cof }F}{\det F}
\ \ \ \ \text{ if }\ \det F>0,
\end{align}
while the case $\det F\le0$ is considered nonphysical. 
In \eq{M-pull-back}, $\bbM:\varOmega\times\R\to\R_{\rm sym}^{d\times d}$ is the
diffusant mobility (depending possibly also on $x\In\varOmega$)
as a material property while ``${\rm Cof}$'' stands for the cofactor matrix.
In literature, this formula is often used in
the isotropic case,
cf.\ e.g.\ \cite[Formula (67)]{DuSoFi10TSMF} or
\cite[Formula (3.19)]{GovSim93CSD2}. For the anisotropic case,
cf.\ \cite{KruRou19MMCM,RouTom13TMPE}.
In fact, \eq{M-pull-back} can be expressed in terms of the right
Cauchy-Green strain $C=F^\top F$ rather than of $F$ itself, which
grants the frame-indifference of this model. The mathematically
interesting attribute of the model \eq{M-pull-back} is that
$\det(\nabla y)$ is (under suitable data qualification) well kept away zero,
similarly as it was already needed for the static problem
because of the capillarity in the actual configuration, and which is
now needed also to \eq{M-pull-back}.


We will use the notation $L^p(0,T;X)$ for the Bochner space of Bochner
measurable functions $[0,T]\to X$ whose norm is in $L^p(0,T)$, 
and $H^1(I;X)$ for functions $[0,T]\to X$ whose distributional derivative
is in $L^2(0,T;X)$. Furthermore, we will not use the Dirichlet condition and 
use the notation $Q=[0,T]\times\varOmega$ and $\Sigma=[0,T]\times\varGamma$.

The departing point is the Hamilton variation principle
adapted for nonconservative systems 
(cf.~also Bedford \cite{Bedf85HPCM}), which says that the integral 
\begin{align}\label{Hamilton-principle}
  \int_0^T\!\!{\mathscr T}(\DT q)-\calE(q(t))+\langle\mathfrak F(t),q(t)\rangle
  \,\d t\ \ \text{ is stationary}
\end{align} 
with $q=(y,\zeta)$ being the state of the system, $\calE$ being the stored
energy and $\mathfrak F(t)=F(t)-\mathcal R'(\DT q)$ being a nonconservative 
force.
This yields the following weak formulation
when one substitutes the concrete functionals:

\begin{definition}[Weak solution]\label{def1}
  The triple $(y,\zeta,\mu)$ with
  $y\in L^\infty(I;H^3(\varOmega;\R^d))\,\cap\,W^{1,\infty}(I;L^2(\varOmega;\R^d))$,
  $\zeta\in L^\infty(I;H^1(\varOmega))$ and $\mu\in L^2(Q)$
  is called a weak solution to the initial-boundary-value
 problem \eq{CH-evol} below if
\begin{subequations}\begin{align}
  &
  \nonumber
\int_Q\big(\pl_F^{}\varphi(\nabla y,\zeta)
+\sigma_\text{\sc k}(\nabla y,\nabla\zeta)\big)\Colon\nabla v
-\varrho\DT{y}\Cdot\DT{v}
+
\nabla^3y\Vdots{\mathbb H}\Vdots\nabla^3v
\,\d x\d t
\\[-.5em]&\hspace{8em}
=\int_Q f\Cdot v\,\d x\d t+\int_{\Sigma}g
\Cdot v\,\d S\d t+\int_\varOmega\varrho v_0\Cdot v(0,\cdot)\,\d x
\label{weak-sln-nonsimple}
\intertext{holds for all 
$v\In L^2(I;H^2(\varOmega;\R^d))\,\cap\,H^1(I;L^2(\varOmega;\R^d))$
with $v|_{t=T}^{}=0$,
and also the initial condition $y(0,\cdot)=y_0$ is satisfied,
and if}
&\nonumber\int_Q\frac{{\mathbb M}(\zeta){\rm Cof}\nabla y}
          {\det\nabla y}\nabla\mu\cdot\big(({\rm Cof}\nabla y)\nabla v\big)-\zeta\DT v\,\d x\d t
          +\int_\Sigma\alpha\mu v\,\d S\d t
\\[-.5em]&\hspace{14em}=\int_\varOmega\zeta_0v(0)\,\d x+\int_\Sigma\alpha\mu_{\rm ext}^{} v\,\d S\d t
\label{weak-sln-zeta}
\intertext{holds for all $v\in H^1(Q)$ with $v|_{t=T}^{}=0$, and}
&\int_Q\kappa(\nabla y)^{\top}\nabla\zeta\cdot(\nabla y)^{\top}\nabla v
+\big(\partial_\zeta^{}\varphi(\nabla y,\zeta)-\mu\big)\, v\,\d x\d t=0
\label{weak-sln-mu}
\end{align}\end{subequations}
for all $v\in L^2(I;H^1(\varOmega))$.
\end{definition}

To see the corresponding initial-boundary-value problem, one is to
apply one by-part integration in time for the inertial term, and here
three-times Green formula over $\varOmega$ and twice surface Green formula
over $\varGamma$.
The resulting boundary-value problem involves rather ``exotic'' hyper-stress
and 3 boundary conditions, namely
\begin{subequations}\label{CH-evol}
\begin{align}\nonumber
&\varrho\DDT y-{\rm div}\big(\pl_F^{}\varphi(\nabla y,\zeta)
+\sigma_\text{\sc k}(\nabla y,\nabla\zeta)
+{\rm div}^2
(\bbH\nabla^3 y)\big)
=f
\label{CH-evol-1}\\&\qquad\qquad\qquad
\text{with }\ \sigma_\text{\sc k}(F,\nabla\zeta)=
\kappa F^{-1}{:}(F^{-\top})'{:}(\nabla\zeta\otimes\nabla\zeta)
&&\text{ in }\ Q,
\\&
\DT\zeta+{\rm div}(\mathfrak{M}(\nabla\bmy,\zeta)\nabla\mu)=0\ \ 
\text{ with }\ 
{\mathfrak M}(\bmF,\zeta)
=\frac{({\rm Cof}\bmF)^\top
  {\mathbb M}(\zeta){\rm Cof}\bmF}
{\det\bmF}&&\text{ in }\ Q,
\label{CH-evol-2}\\&\qquad\qquad\quad
\text{and with }\mu=\pl_z^{}\varphi(\nabla y,\zeta)
-
{\rm div}\big(\kappa(\nabla y)^{-1}(\nabla y)^{-\top}\nabla\zeta\big)
&&\text{ in }\ Q,
\label{CH-evol-3}
\displaybreak\\\nonumber
&
(\pl_F^{}\varphi(\nabla y,\zeta)
+\sigma_\text{\sc k}(\nabla y,\nabla\zeta))\vec{n}
-\divS\big((\divS\vec{n})({\mathbb H}\nabla^3y)\Colon(\vec{n}
\otimes\vec{n})\big)
\\[-.1em]&\nonumber\qquad\qquad\quad\
+\divSS\big((\bbH\nabla^3y)\Cdot\vec{n}\big)
-(\divS\vec{n})\divS\big(({\mathbb H}\nabla^3y)\Cdot\vec{n}\big)\Cdot\vec{n}
\\[-.1em]&\qquad\qquad\qquad\qquad\quad
+{\rm div}^2({\mathbb H}\nabla^3y)\Cdot\vec{n}
+\divS\big({\rm div}({\mathbb H}\nabla^3y)\Cdot\vec{n}\big)
=g
&&\text{ on }\ \Sigma,
\label{CH-evol-4}\\[-.5em]&
(\bbH\nabla^3y)\Vdots(\vec{n}\otimes\vec{n}\otimes\vec{n})=0
\ \ \text{ and }\ \ {\rm div}(\bbH\nabla^3y)\Vdots(\vec{n}\otimes\vec{n})=0
,\ \ 
&&\text{ on }\ \Sigma,
\label{CH-evol-5}\\&
\mathfrak{M}(\nabla\bmy,\zeta)(\nabla\zeta\otimes\vec{n})
+\alpha\mu=\alpha\mu_{\rm ext}\ \ &&\text{ on }\ \Sigma,
\label{CH-evol-6}
\\&\kappa(\nabla y)^{-1}(\nabla y)^{-\top}\nabla\zeta\Cdot\vec{n}=0
&&\text{ on }\ \Sigma\,,
\label{CH-evol-7}
\end{align}\end{subequations}
where $\divS={\rm tr}(\nablaS)$ with ${\rm tr}(\cdot)$ being the trace of a
$(d{-}1){\times}(d{-}1)$-matrix, denotes the $(d{-}1)$-dimensional
surface divergence and
$\nablaS v=\nabla v-\frac{\partial v}{\partial\vec{n}}\vec{n}$ 
being the surface gradient of $v$.
Here we use, in addition what would come from \eq{Hamilton-principle}
with \eq{functionals},
also a nonhomogeneous boundary condition for the diffusion, involving
an external chemical potential $\mu_{\rm ext}$.
The variable $\mu$ from \eq{CH-evol-3} is called a \emph{chemical potential}
and $\mathfrak{M}(\nabla\bmy,\zeta)\nabla\mu$ in \eq{CH-evol-2} is the
Fick law for the flux of the diffusant.

The system \eq{CH-evol} deserves some comments. First, the diffusion equation
(\ref{CH-evol}b,c) considered with the 
Robin boundary conditions \eq{CH-evol-6} 
$\mathfrak{M}(\nabla y,\zeta)\nabla\mu\Cdot\vec{n}+\alpha\mu=\mu_{\rm ext}^{}$ on
$\varGamma$ can be rewritten in the form
\begin{align}\label{poro-evol-flow}
  \Delta_{\mathfrak{M}(\nabla y,\zeta)}^{-1}(\DT\zeta,\mu_{\rm ext}^{})=\mu=\pl_z^{}\varphi(\nabla y,\zeta)
-{\rm div}\big(\kappa(\nabla y)^{-1}(\nabla y)^{-\top}\nabla\zeta\big)
\end{align}
with $\Delta_{\mathfrak{M}(\nabla y,\zeta)}^{-1}$ as in \eq{diss-pot}. 
In view of (\ref{functionals}a,c), this is exactly
$\pl_{\DT\zeta}\mathcal R(y,\zeta;\DT\zeta)+\pl_\zeta^{}\calE(y,\zeta)=0$, which
is what results from the Hamilton variational principle
\eq{Hamilton-principle} as far as $\zeta$-component concerns.

Further, the boundary conditions (\ref{CH-evol}d,e) for the mechanical
equilibrium \eq{CH-evol-1} are quite technical
because of the $\bbH$-term.
It is to be treated, at each time instant $t$ (not explicitly denoted),
first by applying three times Green formula
\begin{align}\nonumber
&
  \int_\varOmega\nabla^3y\Vdots{\mathbb H}\Vdots\nabla^3v\,\d x
  =\int_\varGamma\big(({\mathbb H}\nabla^3y)\Cdot\vec{n}\big)\Vdots\nabla^2v\,\d S
  -\int_\varOmega{\rm div}({\mathbb H}\nabla^3y)\Vdots\nabla^2v\,\d x
  \\&\nonumber
  =\int_\varGamma\big(({\mathbb H}\nabla^3y)\Cdot\vec{n}\big)\Vdots
  \big(\partial_{\vec{n}}^2v{+}\nablaSS v\big)\,\d S
  -\int_\varOmega{\rm div}({\mathbb H}\nabla^3y)\Vdots\nabla^2v\,\d x 
  \\&\nonumber
  =\int_\varGamma\big(({\mathbb H}\nabla^3y)\Cdot\vec{n}\big)\Vdots
  \big(\partial_{\vec{n}}^2v{+}\nablaSS v\big)
  \\[-.5em]&\nonumber\qquad\qquad
  -\big({\rm div}({\mathbb H}\nabla^3y)\Cdot\vec{n}\big)\Colon
  (\partial_{\vec{n}}^{}v{+}\nablaS v)\,\d S
  +\int_\varOmega{\rm div}^2({\mathbb H}\nabla^3y)\Colon\nabla v\,\d x 
  \\[-.3em]&\nonumber
  =\int_\varGamma\big(({\mathbb H}\nabla^3y)\Cdot\vec{n}\big)\Vdots
  \big(\partial_{\vec{n}}^2v{+}\nablaSS v\big)
  +\big({\rm div}^2({\mathbb H}\nabla^3y)\Cdot\vec{n}\big)\Cdot v
  \\[-.5em]&\qquad\qquad
  -\big({\rm div}({\mathbb H}\nabla^3y)\Cdot\vec{n}\big)\Colon
  (\partial_{\vec{n}}^{}v{+}\nablaS v)\,\d S
  -\int_\varOmega{\rm div}^3({\mathbb H}\nabla^3y)\Cdot v\,\d x\,,
\label{H-calcus}\end{align}
where we used the decomposition of $\nabla v$ on $\varGamma$ into the normal and the
tangential part $\partial_{\vec{n}}^{}v+\nablaS v$, and in particular also
$$\nabla^2v=
(\partial_{\vec{n}}^{}+\nablaS)(\partial_{\vec{n}}^{}v+\nablaS v)
=\partial_{\vec{n}}^2v+\nablaSS v\,,
$$
where we use also the orthogonality of $\nablaS v$ and $\partial_{\vec{n}}^{}v$.
We further apply four times the surface Green formula on the boundary term;
more specifically, we apply 
\begin{align}\label{eq:Green-surf}
  \int_{\varGamma}A{:}\nablaS v\,\d S
  =\int_\varGamma((\divS\vec{n})A\vec{n}-\divS A)\cdot v\,\d S\,, 
\end{align} 
which holds for a smooth field $A\in C^1(\varGamma;\R^{d\times d})$ and
$v\in C^1(\varGamma;\R^d)$ that; cf.\ see \cite[Formula (34)]{FriGur06TBBC}.
By this way, we can write
\begin{align}
&\nonumber
  \!\!\!\int_\varGamma\big(({\mathbb H}\nabla^3y)\Cdot\vec{n}\big)\Vdots
  \big(\partial_{\vec{n}}^2v{+}\nablaSS v\big)
  +\big({\rm div}^2({\mathbb H}\nabla^3y)\Cdot\vec{n}\big)\Cdot v
  -\big({\rm div}({\mathbb H}\nabla^3y)\Cdot\vec{n}\big)\Colon
  (\partial_{\vec{n}}^{}v{+}\nablaS v)\,\d S
\\[-.3em]&\nonumber
\!\!\!=\int_\varGamma
\big(({\mathbb H}\nabla^3y)\Cdot\vec{n}\big)\Vdots\partial_{\vec{n}}^2v
-\big({\rm div}({\mathbb H}\nabla^3y)\Cdot\vec{n}\big)\Colon\partial_{\vec{n}}^{}v
\\[-.6em]&\nonumber\qquad
+\Big((\divS\vec{n})({\mathbb H}\nabla^3y)\Colon(\vec{n}\otimes\vec{n})\big)
-\divS\big(({\mathbb H}\nabla^3y)\Cdot\vec{n}\big)
\Big)\Colon\nablaS v
\\[-.1em]&\nonumber\qquad
+\Big({\rm div}^2({\mathbb H}\nabla^3y)\Cdot\vec{n}
-(\divS\vec{n}){\rm div}({\mathbb H}\nabla^3y)\Colon(\vec{n}\otimes\vec{n})
+\divS\big({\rm div}({\mathbb H}\nabla^3y)\Cdot\vec{n}\big)
\Big)\Cdot v\,\d S
\\&\nonumber
\!\!\!=\int_\varGamma\big(({\mathbb H}\nabla^3y)\Cdot\vec{n}\big)\Vdots\partial_{\vec{n}}^2v
-\big({\rm div}({\mathbb H}\nabla^3y)\Cdot\vec{n}\big)\Colon\partial_{\vec{n}}^{}v
\\[-.8em]&\nonumber\qquad
+\Big((\divS\vec{n})^2({\mathbb H}\nabla^3y)\Vdots(\vec{n}\otimes\vec{n}\otimes\vec{n})\big)
-\divS\big((\divS\vec{n})({\mathbb H}\nabla^3y)\Colon(\vec{n}
\otimes\vec{n})\big)
\\[-.1em]&\nonumber\qquad
+\divSS\big((\bbH\nabla^3y)\Cdot\vec{n}\big)
-(\divS\vec{n})\divS\big(({\mathbb H}\nabla^3y)\Cdot\vec{n}\big)\Cdot\vec{n}
+{\rm div}^2({\mathbb H}\nabla^3y)\Cdot\vec{n}
\\[-.1em]&\qquad
-(\divS\vec{n}){\rm div}({\mathbb H}\nabla^3y)\Colon(\vec{n}\otimes\vec{n})
+\divS\big({\rm div}({\mathbb H}\nabla^3y)\Cdot\vec{n}\big)
\Big)\Cdot v\,\d S.
\label{H-calcus+}\end{align}
Substituting \eq{H-calcus+} into \eq{H-calcus} and into
\eq{weak-sln-nonsimple}, and taking $v$ arbitrarily
with compact support, then with arbitrary traces but with normal derivatives
zero, and then with $\partial_{\vec{n}}^2v=0$, and eventually entirely
arbitrarily, we obtain subsequently \eq{CH-evol-1} and (\ref{CH-evol}d,e).
Notably, the conditions \eq{CH-evol-5} have been reflected also 
in \eq{CH-evol-4} to simplify it in contrast what can be seen
from the last seven terms in \eq{H-calcus+}.

It is important that this gradient theory in the actual configuration
has led to a specific contribution $\sigma_\text{\sc k}$ to 
the stress tensor. Such stresses are needed, in particular,  to 
balance energy and are known in incompressible-fluid mechanics under the name 
\emph{Korteweg stresses} \cite{Kort01FPEM}. Evaluating $(F^{-\top})'$ and
eliminating $F^{-1}={\rm Cof}F^\top/\det F$, this stress can be expressed more
specifically as 
\begin{align}
  \sigma_\text{\sc k}(F,\nabla\zeta)=\kappa
  \frac{{\rm Cof}F^\top\!\!}{\det F}\bigg(\frac{{\rm Cof}'F}{\det F}-
\frac{{\rm Cof}F\otimes{\rm Cof}F}{(\det F)^2}\bigg)
\Colon(\nabla\zeta\otimes\nabla\zeta)\,.
\end{align}
Mathematically, this stress brings an additional difficulty in comparison with
the usual concepts of gradients in reference configuration, because $\nabla\zeta$
occurs nonlinearly and we need the strong convergence of an approximation of
$\nabla\zeta$.

\begin{proposition}[Existence of weak solutions to the poro-elasto-dynamics.]
  Let $d=2,3$ and \eq{ass-growrt-phi-dam} hold for $p>2d/(d{-}2)$,
  and let $\mathbb M$ be a bounded Carath\'eodory mapping with values
  uniformly positive definite. Moreover, let $f\In L^1(I;L^2(\varOmega;\R^d))$, 
  $g\In W^{1,1}(I;L^1(\varGamma;\R^d))$, $\mu_{\rm ext}\In L^2(\Sigma)$,
  $y_0\in H^3(\varOmega;\R^d)$, $\inf_\varOmega^{}\det\nabla y_0>0$,
  $v_0\in L^2\varOmega;\R^d)$, and $\zeta_0\in H^1(\varOmega)$.
  Then there exists a weak solution $(y,\zeta,\mu)$ to the
  initial-boundary-value problem \eq{CH-evol} according Definition~\ref{def1}
such that
 $\varrho\DDT y\in L^1(I;H^3(\varOmega;\R^d)^*)$ and
 $\DT\zeta\in L^2(I;H^1(\varOmega)^*)$.
\end{proposition}

  \noindent{\it Proof.}
  We first construct the conformal Galerkin approximation of \eq{CH-evol}.
  This means the finite-dimensional subspaces for \eq{CH-evol-1}
  are contained in $H^3(\varOmega;\R^d)$, while for \eq{CH-evol-2} and
  \eq{CH-evol-3} they are contained in $H^1(\varOmega)$. Let us denote the
  solution obtained by this way as $(y_k,\zeta_k,\mu_k)$ with $k\in\N$ denoting
  the indexing of the mentioned finite-dimensional subspaces. Existence
  of such approximate solution is by usual continuation argument, based on
  the uniform a-priori estimates below.

  It is important to take these finite-dimensional subspaces for \eq{CH-evol-2}
  and for \eq{CH-evol-3} (written for the approximate solution) the same in
  order to allow for a cross-test of \eq{CH-evol-2} by $\mu_k$ and
  \eq{CH-evol-3} by $\DT\zeta_k$.

  Together with the test of \eq{CH-evol-1} by $\DT y_k$ and using the boundary
  conditions (\ref{CH-evol}d--g), we obtain the discrete energy balance
\begin{align}\nonumber
 &
\frac{\d}{\d t}\int_\varOmega\frac\varrho2|\DT y_k|^2+\varphi(\nabla\bmy_k,\zeta_k)
+\frac{\kappa}2\Big|
\frac{{\rm Cof}\nabla\bmy_k}{\det\nabla\bmy_k}\nabla\zeta_k\Big|^2
  +\frac12\bbH\nabla^3\bmy_k\Vdots\nabla^3\bmy_k\,\d x
  \\&\qquad\qquad\nonumber +\int_\varOmega 
\frac{({\rm Cof}\nabla\bmy_k)^\top{\mathbb M}(\zeta_k){\rm Cof}\nabla\bmy_k)}
     {\det\nabla\bmy_k}\nabla\mu_k\Cdot\nabla\mu_k\,\d x
     +\int_\varSigma\alpha\mu_k^2\,\d S
     \\[-.3em]&\qquad\qquad\qquad\qquad\qquad\qquad\quad
     =\int_\varOmega f\Cdot\DT y_k\,\d x+\int_\Sigma g\Cdot\DT y_k+
\alpha\mu_{\rm ext}\mu_k\,\d S\,.
\label{energy-k}\end{align}
 Here we have enjoyed cancellation of the terms $\pm\mu_k\DT c_k$ and have used 
the calculus $\plF\varphi(\nabla y_k,\zeta_k)\Colon\nabla\DT y_k
+\pl_z^{}\varphi(\nabla y_k,c_k)\DT \zeta_k
=\frac{\pl}{\pl t}\varphi(\nabla y_k,\zeta_k)$.

We integrate \eq{energy-k} over time interval $[0,t]$ and apply
by-part integration in time on the term $g\Cdot\DT y_k$ because
$\DT y_k$ does not have well estimated traces on $\varGamma$. Then 
we apply the H\"older and the Gronwall inequalities.
By the Healey-Kr\"omer theorem \cite{HeaKro09IWSS}
holding, in fact, on each level sets and being here
used with the compact embedding $H^3(\varOmega;\R^d)\subset W^{2,p}(\varOmega;\R^d)$,
we have
\begin{align}
\forall (t,x)\In\bar Q:\qquad\det\nabla y_k(t,x)\ge\varepsilon
\end{align}
for some positive 
$\varepsilon\le\min_{x\in\bar\varOmega}\det\nabla y_0(x)$
It is important that this holds by successive-continuation argument on the
Galerkin level, and thus $\nabla y_k$ is valued in the definition domain
of $\varphi$ and the singularity of $\varphi$ is not seen, and therefore
the Lavrentiev phenomenon is excluded. 
Altogether, by this way, we obtain the a-priori estimates
\begin{subequations}\label{IBVP-diffusion-large-est}
\begin{align}\label{IBVP-diffusion-large-est1}
&\|y_k\|_{L^\infty(I;H^3(\varOmega;\R^d))\,\cap\,W^{1,\infty}(I;L^2(\varOmega;\R^d))}^{}\le K\ \ \text{ and }\ \ \Big\|\frac1{\det\nabla y_k}\Big\|_{L^\infty(Q)}\le K,
\\&\label{IBVP-diffusion-large-est3}
\bigg\|\frac{{\rm Cof}\nabla y_k}{\sqrt{\det\nabla y_k}}\nabla\mu_k\bigg\|_{L^2(Q;\R^d)}^{}\!\le K,
\\&\label{IBVP-diffusion-large-est4}
\bigg\|
\frac{{\rm Cof}\nabla\bmy_k}{\det\nabla\bmy_k}\nabla\zeta_k
\bigg\|_{L^\infty(I;L^2(\varOmega;\R^d))}\!\le K.
\end{align}\end{subequations}
From \eq{IBVP-diffusion-large-est1}, we have the bound
$\nabla y_k\In L^\infty(Q;\R^{d\times d})
$ so that, realizing that $({\rm Cof}\nabla y_k)^{-1}=
(\nabla y_k)^\top/\det\nabla y_k$, from \eq{IBVP-diffusion-large-est3}
we have
\begin{align}\nonumber
&\|\nabla\mu_k\|_{L^2(Q
;\R^d))}^{}=
\bigg\|(\nabla y_k)^\top\frac{{\rm Cof}\nabla y_k}{\det\nabla y_k}\nabla\mu_k\bigg\|_{L^2(
Q;\R^d))}^{}
\\&\qquad\qquad\le\|\nabla y_k\|_{L^\infty(Q;\R^{d\times d})}^{}
\bigg\|\frac{1}{\sqrt{\det\nabla y_k}}\bigg\|_{L^\infty(Q)}
\bigg\|\frac{({\rm Cof}\nabla y_k)\nabla\mu_k}{\sqrt{\det\nabla y_k}}\bigg\|_{L^2(Q;\R^{d\times d})}^{}.
\end{align}
Then we use \eq{IBVP-diffusion-large-est3}
to obtain the bound of $\nabla\mu_k$ in $L^2(Q;\R^d)$.
By the Poincar\'e inequality based on
the Robin boundary condition we obtain the bound of
\begin{align}\label{est-muk}
 \|\mu_k\|_{L^2(I;H^1(\varOmega))}^{}\le K\,.
\end{align}
Similarly, we can estimate
\begin{align}\nonumber
\|\nabla\zeta_k\|_{L^\infty(I;L^2(\varOmega;\R^d))}^{}&=
  \bigg\|(\nabla y_k)^\top\frac{{\rm Cof}\nabla y_k}{\det\nabla y_k}\nabla\zeta_k\bigg\|_{L^\infty(I;L^2(\varOmega;\R^d))}^{}
\\&\le
\|\nabla y_k\|_{L^\infty(Q;\R^{d\times d})}^{}\bigg\|\frac{{\rm Cof}\nabla y_k}{\det\nabla y_k}\nabla\zeta_k\bigg\|_{L^\infty(I;L^2(\varOmega;\R^d))}^{}\,,
\end{align}
from which we obtain the bound of $\nabla\zeta_k$ in
$L^\infty(I;L^2(\varOmega;\R^d))$ by using \eq{IBVP-diffusion-large-est3}.
Then, by the coercivity of $\varphi(F,\cdot)$, cf.\ \eq{growth}, also the
estimate
\begin{align}
\label{IBVP-diffusion-large-est2}
\|\zeta_k\|_{L^\infty(I;H^1(\varOmega))}^{}\le K.
\end{align}

Then we select a weakly* convergent subsequence in the topologies indicated in
\eq{IBVP-diffusion-large-est1},
\eq{est-muk}, and \eq{IBVP-diffusion-large-est2}. Moreover, by comparison,
from the equation \eq{CH-evol-2} in its Galerkin
approximation and from \eq{IBVP-diffusion-large-est3},
we can also see that (a Hahn-Banach extension of) $\DT\zeta_k$ is bounded in 
$L^2(I;H^1(\varOmega)^*)$. Then one can used the Aubin-Lions lemma
to get strong convergence
both for 
$$
\nabla y_k\to\nabla y\ \ \text{ in }\ C(\bar Q;\R^{d\times d})
\ \ \text{ and }\ \ \zeta_k\to\zeta\text{ in }\ L^{1/\varepsilon}(I;L^{p-\varepsilon}(\varOmega))
$$
for any
$0<\varepsilon<p-1$ with $p=6$ if $d=3$ or $p<+\infty$ if $d=2$.
The convergence towards the weak 
solution of \eq{CH-evol} is then easy.

A bit peculiar term is the diffusion flux when considering
the ansatz \eq{M-pull-back} and thus the weak formulation 
\eq{weak-sln-zeta}, for which we need to show that
\begin{align}\nonumber
&\int_Q\bigg(\bbM(\zeta_k)\frac{{\rm Cof}\nabla y_k}{\sqrt{\det\nabla y_k}}\nabla\mu_k\bigg)
\Cdot\bigg(\frac{{\rm Cof}\nabla y_k}{\sqrt{\det\nabla y_k}}\nabla v\bigg)\,\d x\d t
\\[-.0em]&\qquad\qquad\qquad\qquad\to
\int_Q\bigg(\bbM(\zeta)\frac{{\rm Cof}\nabla y}{\sqrt{\det\nabla y}}\nabla\mu\bigg)
\Cdot\bigg(\frac{{\rm Cof}\nabla y}{\sqrt{\det\nabla y}}\nabla v\bigg)\,\d x\d t
\end{align}
for any $v\In C^1(\bar Q)$. Here we used that 
\begin{subequations}\begin{align}
&\frac{{\rm Cof}\nabla y_k}{\sqrt{\det\nabla y_k}}\nabla\mu_k\to
\frac{{\rm Cof}\nabla y}{\sqrt{\det\nabla y}}\nabla\mu\!\!\!\!&&\text{weakly in }\ 
L^2(Q;\R^{d}),\ \ \text{ and}
\\& 
\frac{{\rm Cof}\nabla y_k}{\sqrt{\det\nabla y_k}}\to
\frac{{\rm Cof}\nabla y}{\sqrt{\det\nabla y}}&&\text{strongly in }\ 
C(\bar Q;\R^{d\times d})
\end{align}\end{subequations}
because ${\rm Cof}\nabla y_k\to {\rm Cof}\nabla y$ strongly in 
$L^p(Q;\R^{d\times d})$ and $1/\det\nabla y_k\to 1/\det\nabla y$ strongly in 
$L^p(Q)$ for any $1\le p<+\infty$ due to  
the Aubin-Lions theorem together with the latter estimate in 
\eq{IBVP-diffusion-large-est1}.
and that
$\bbM(c_k)(\nabla y_k)^{-1}\to\bbM(c)(\nabla y)^{-1}$ weakly in $L^2(Q)$
thanks to the estimate \eq{IBVP-diffusion-large-est3}. 

  As already mentioned, the limit passage in the Korteweg-like stress
  $\sigma_\text{\sc k}$ needs strong convergence of $\nabla\zeta_k$ in
  $L^2(Q;\R^d)$. To this goal, 
we use the uniform (with respect to $y$) strong monotonicity of the mapping
$\zeta\mapsto-{\rm div}\big(\kappa(\nabla y)^{-1}(\nabla y)^{-\top}\nabla\zeta\big)$.
Taking $\widetilde\zeta_k$ an approximation of $\zeta$ valued in the
respective finite-dimensional spaces used for the Galerkin approximation
and converging to $\zeta$ strongly, we can test \eq{CH-evol-3} in its Galerkin
approximation by $\zeta_k{-}\widetilde\zeta_k$ and use it in the estimate
\begin{align}\nonumber
&\!\!\limsup_{k\to\infty}\int_Q\kappa(\nabla y_k)^{-1}(\nabla y_k)^{-\top}
\nabla(\zeta_k{-}\widetilde\zeta_k)\Cdot\nabla(\zeta_k{-}\widetilde\zeta_k)
\,\d x\d t
\\&\nonumber\qquad=\lim_{k\to\infty}\int_Q\big(\pl_z\varphi(\nabla y_k,\zeta_k)
+\mu_k\big)(\widetilde\zeta_k{-}\zeta_k)
\\[-.3em]&\qquad\qquad\qquad\
-\kappa(\nabla y_k)^{-1}(\nabla y_k)^{-\top}\nabla\widetilde\zeta_k
\Cdot\nabla(\zeta_k{-}\widetilde\zeta_k)\,\d x\d t=0
\label{nabla-zeta-strong}\end{align}
because $\pl_z\varphi(\nabla y_k,\zeta_k)+\mu_k$ is
bounded in $L^2(Q)$ while $\widetilde\zeta_k{-}\zeta_k\to0$ strongly in
$L^2(Q)$ by the Aubin-Lions compactness theorem and because 
$\kappa(\nabla y_k)^{-1}(\nabla y_k)^{-\top}\nabla\widetilde\zeta_k$ converges
strongly in $L^2(Q;\R^d)$ while $\nabla(\zeta_k{-}\widetilde\zeta_k)\to0$
weakly in $L^2(Q;\R^d)$. As $\kappa(\nabla y_k)^{-1}(\nabla y_k)^{-\top}$ is
uniformly positive definite, we thus obtain that
$\nabla(\zeta_k{-}\widetilde\zeta_k)\to0$ strongly in $L^2(Q;\R^d)$,
and thus $\nabla\zeta_k\to\nabla\zeta$ strongly in $L^2(Q;\R^d)$.

Then we have the convergence in the Korteweg-like stress even strongly in
$L^p(I;L^1(\varOmega;\R^{d\times d}))$ for any $1\le p<+\infty$. The limit passage
in the force equilibrium towards (\ref{CH-evol}a,d,e) formulated weakly in
\eq{weak-sln-nonsimple} is the straightforward.
$\hfill\Box$

\section{Concluding remarks}\label{sect-remarks}
We close the paper with a brief outlook to some modifications and other applications
and models which can be analysed quite analogously.

\begin{remark}[\emph{Allen-Cahn modification: damage or phase-transformation
      models.}]\upshape Replacing the quadratic (in terms of rate) nonlocal dissipation potential
  by a non-quadratic nonsmooth (at zero-rate) local 
  dissipation potential of the type
\begin{align*}
  R(\DT\zeta):=\int_\varOmega r(\DT\zeta)\,\d x\quad\text{ with }\
  r:\R\to[0,+\infty]\ \text{ convex},
\end{align*}
we would obtain a diffusionless model of Allen-Cahn type \cite{AllCah72GSSO}.
The equation \eq{CH-evol-2} is then simplified for $\pl r(\DT\zeta)+\mu=0$ and
\eq{CH-evol-6} is omitted, while the Korteweg-like contribution $\sigma_{\sc k}$
induced by the actual-configuration gradient of $\zeta$ to the stress tensor
remains in \eq{CH-evol-1}. This may describe a damage model \cite{WPMB14GELD}
or a martensitic phase transformation \cite{Levi14PFAM}.
The mentioned nonsmoothness of $r(\cdot)$ at $\DT\zeta=0$ then models
activation phenomena and, in the case of reversible phase transformation,
hysteresis behaviour. For the analysis, we refer to
\cite[Sect.\,9.5.1]{KruRou19MMCM}.
\end{remark}

\begin{remark}[\emph{Dispersion of elastic waves.}]\upshape
  The concept of nonsimple materials allow for introducing a dispersion
  of elastic waves, as well known for linear models at small strains.
  Typically, involving higher gradients in a positive-definite way,
  one gets anomalous dispersion, i.e.\ higher-frequency waves propagate
  faster than lower-frequency ones.
When one combines the concept of 3rd-grade (as here in Sect.~\ref{sect-dynam})
with the 2nd-grade (as in Sect.~\ref{sect-static}) materials, we obtain
a bigger freedom. In particular, a combination of normal and anomalous dispersion
can be obtained when the second-order deformation gradient is involved
in a negative-definite way, cf. also \cite[Remark~6.3.6]{KruRou19MMCM}
for a 1-dimensional linear model.
\end{remark}


\begin{remark}[\emph{Gradient plasticity.}]
\label{rem-plast-grad-actual}\upshape
Another model where gradient can be considered in the actual deforming
configuration is plasticity. At large strains, it is always analytically
necessary to involve gradient of plastic strain $\PP$ into the stored energy,
which is then considered as
\begin{align}\label{plast-large-stored+}
\calE(\nabla y,\PP)=\int_\varOmega\varphi_\text{\sc el}((\nabla y)\PP^{-1})
+\varphi_\text{\sc hd}(\PP)+\frac12\nabla^3y\Vdots\bbH\Vdots\nabla^3y
+\frac\kappa p\big|(\nabla y)^{-\top}\nabla\PP\big|^p.
\end{align}
When considering still the kinetic energy \eq{T} and the dissipation
potential $\mathcal{R}(\PP;\DT\PP)=\int_\varOmega\zeta(\DT\PP\PP^{-1})\,\d x$
for some convex $\zeta:\R^{d\times d}\to\R$,
the evolution system arising by the Hamilton variational principle
extends as
\begin{subequations}\label{plast-plast-large}
\begin{align}\label{plast-heat-large-1}
&\varrho\DDT y-{\rm div}\,\SSel=f(y)\ \ \ \text{ with }\ \
&&\text{ in }\ Q,
\\&\nonumber
\partial_P^{}{\zeta}\big(\DT\PP\PP^{-1}\big)+\,\SSin\PP^\top
\ni0\ \ \ &&\text{ in }\ Q,
\end{align}\end{subequations}
with the elastic stress 
$\SSel=
\pl_{\nabla y}^{}\calE(\nabla y,\PP)$ and
an inelastic driving stress 
$\SSin=\pl_\PP^{}\calE(\nabla y,\PP)$. In view of \eq{plast-large-stored+}, 
we can specify
\begin{subequations}\begin{align}\nonumber
    &\SSel=\varphi_\text{\sc el}'(\nabla y\PP^{-1})\PP^{-\top}+{\rm div}^2(\bbH\nabla^3 y)
   \\&\qquad +\kappa|(\nabla y)^{-\top}\nabla\PP\big|^{p-2}\nabla\PP^{\top}(\nabla y)^{-1}(\nabla y)^{-\top}\nabla\PP\quad\text{ and }
     \label{Korteweg-plast}
     \\&\nonumber
\SSin=(\nabla y)^\top\varphi_\text{\sc el}'(\nabla y\PP^{-1})\Colon(\PP^{-1})'
+\varphi_\text{\sc hd}'(\PP)
\\&\qquad
-{\rm div}\big(\kappa|(\nabla y)^{-\top}\nabla\PP\big|^{p-2}(\nabla y)^{-1}
(\nabla y)^{-\top}\nabla\PP\big)\end{align}\end{subequations}
The last term in \eq{Korteweg-plast}
is a \emph{Korteweg}-like
\emph{stress}\index{stress!Korteweg}\index{Korteweg stress} and,
because of it, now the strong convergence in $\nabla\PP$ is needed for
the analysis. This can be done similarly as \eqref{nabla-zeta-strong}, now 
based on the uniform (with respect to $\nabla y$) strong monotonicity
of the mapping $\PP\mapsto
-{\rm div}(\kappa|(\nabla y)^{-\top}\nabla\PP\big|^{p-2}(\nabla y)^{-1}
(\nabla y)^{-\top}\nabla\PP)$; cf.\ \cite[Sect.\,9.4.2]{KruRou19MMCM}.
A combination of the Cahn-Hilliard models with plasticity can also be
considered, like \cite{Anan12CHTT,RouSte18TEPR}.
\end{remark}


\begin{remark}[Open problems.]\upshape
The gradient of the deformation gradient in
\eq{stored-static} and \eq{stored-dynam} is considered in the
reference configuration while the concentration gradient is in
the actual deformed configuration. This is a certain conceptual
discrepancy. Yet, considering the non-simple materials in the
actual configuration brings additional terms and serious additional
difficulties. \COMMENT{DETAILS? static MAYBE possible ??}
\end{remark}

%

\begin{thebibliography}{10}

\bibitem{AgiLaz09CBLL}
E.~Agiasofitou and M.~Lazar.
\newblock Conservation and balance laws in linear elasticity of grade three.
\newblock {\em J. Elasticity}, 94:69--85, 2009.

\bibitem{AllCah72GSSO}
S.~Allen and J.~Cahn.
\newblock Ground state structures in ordered binary alloys with second neighbor
  interactions.
\newblock {\em Acta Metall.}, 20:423--433, 1972.

\bibitem{Anan12CHTT}
L.~Anand.
\newblock A {C}ahn-{H}illiard-type theory for species diffusion coupled with
  large elastic-plastic deformations.
\newblock {\em J. Mech. Phys. Solids}, 60:1983--2002, 2012.

\bibitem{ArSaRa16SACC}
P.~Areias, E.~Samaniego, and T.~Rabczuk.
\newblock A staggered approach for the coupling of {C}ahn-{H}illiard type
  diffusion and finite strain elasticity.
\newblock {\em Comput. Mech.}, 57:339--351, 2016.

\bibitem{Bedf85HPCM}
A.~Bedford.
\newblock {\em Hamilton's Principle in Continuum Mechanics}.
\newblock Pitman, Boston, 1985.

\bibitem{BCDGSS02MPSBA}
E.~Bonetti, P.~Colli, W.~Dreyer, G.~Gilardi, G.~Schimperna, and J.~Sprekels.
\newblock On a model for phase separation in binary alloys drivenby mechanical
  effects.
\newblock {\em Physica D}, 165:48--65, 2002.

\bibitem{CahHil58FEUS}
J.~Cahn and J.~Hilliard.
\newblock Free energy of a uniform system {I}., {I}nterfacial free energy.
\newblock {\em J. Chem. Phys.}, 28:258--267, 1958.

\bibitem{CiaNec87ISCN}
P.~Ciarlet and J.~Ne{\v{c}}as.
\newblock Injectivity and self-contact in nonlinear elasticity.
\newblock {\em Arch. Ration. Mech. Anal.}, 97:171--188, 1987.

\bibitem{DalMie15CECM}
H.~Dal and C.~Miehe.
\newblock Computational electro-chemo-mechanics of lithium-ion battery
  electrodes at finite strains.
\newblock {\em Comput. Mech.}, 55:303--325, 2015.

\bibitem{DLReAn14CHTP}
C.~{Di Leo}, E.~Rejovitzky, and L.~Anand.
\newblock A {C}ahn-{H}illiard-type phase-field theory for species diffusion
  coupled with large elastic deformations: Application to phase-separating
  {L}i-ion electrode materials.
\newblock {\em J. Mech. Phys. Solids}, 70:1--29, 2014.

\bibitem{DuSoFi10TSMF}
F.~Duda, A.~Souza, and E.~Fried.
\newblock A theory for species migration in a finitely strained solid with
  application to polymer network swelling.
\newblock {\em J. Mech. Phys. Solids}, 58:515--529, 2010.

\bibitem{FriGur06TBBC}
E.~Fried and M.~Gurtin.
\newblock Tractions, balances, and boundary conditions for nonsimple materials
  with application to liquid flow at small-lenght scales.
\newblock {\em Arch. Ration. Mech. Anal.}, 182:513--554, 2006.

\bibitem{Garc03CHSE}
H.~Garcke.
\newblock On {C}ahn-{H}illiard system with elasticity.
\newblock {\em Proc. Roy. Soc. Edinburgh}, 133A:307--331, 2003.

\bibitem{GovSim93CSD2}
S.~Govindjee and J.~Simo.
\newblock Coupled stress-diffusion: case {II}.
\newblock {\em J. Mech. Phys. Solids}, 41:863--887, 1993.

\bibitem{HeaKro09IWSS}
T.~Healey and S.~Kr\"omer.
\newblock Injective weak solutions in second-gradient nonlinear elasticity.
\newblock {\em ESAIM: Control, Optim. \& Cal. Var.}, 15:863--871, 2009.

\bibitem{HeiKra14PSCD}
C.~Heinemann and C.~Kraus.
\newblock {\em Phase Separation Coupled with Damage Processes}.
\newblock Springer, Wiesbaden, 2014.

\bibitem{HGOD17FPLS}
C.~Hesch, A.~Gil, R.~Ortigosa, M.~Dittmann, C.~Bilgen, P.~Betsch, M.~Franke,
  A.~Janz, and K.~Weinberg.
\newblock A framework for polyconvex large strain phase-field methods to
  fracture.
\newblock {\em Comput. Methods Appl. Mech. Engrg.}, 317:649--683, 2017.

\bibitem{HonWan13PFMS}
W.~Hong and X.~Wang.
\newblock A phase-field model for systems with coupled large deformation and
  mass transport.
\newblock {\em J. Mech. Phys. Solids}, 61:1281--1294, 2013.

\bibitem{Kort01FPEM}
D.~Korteweg.
\newblock Sur la forme que prennent les \'equations du mouvement des fuides si
  l\'on tient compte des forces capillaires caus\'ees par des variations de
  densit\'e consid\'erables mais continues et sur la th\'eorie de la
  capillarit\'e dans l'hypoth\`ese d'une variation continue de la densit\'e.
\newblock {\em Arch. N\'eerl. Sci. Exactes Nat.}, 6:1--24, 1901.

\bibitem{KruRou19MMCM}
M.~Kru\v{z}\'{\i}k and T.~Roub\'{\i}\v{c}ek.
\newblock {\em Mathematical Methods in Continuum Mechanics of Solids}.
\newblock Springer, Switzerland, 2019.

\bibitem{LarCah82ESSD}
F.~Larch\'{e} and J.~Cahn.
\newblock The effect of self–stress on diffusion in solids.
\newblock {\em Acta Metall.}, 30:1835--1845, 1982.

\bibitem{Levi14PFAM}
V.~Levitas.
\newblock Phase field approach to martensitic phase transformations with large
  strains and interface stresses.
\newblock {\em J. Mech. Phys. Solids}, 70:154--189, 2014.

\bibitem{MiMaUl14FNEM}
C.~Miehe, S.~Mauthe, and H.~Ulmer.
\newblock Formulation and numerical exploitation of mixed variational
  principles for coupled problems of {C}ahn-{H}illiard-type and standard
  diffusion in elastic solids.
\newblock {\em Int. J. Numer. Meth. Engng.}, 99:737--762, 2014.

\bibitem{Miel11GSRD}
A.~Mielke.
\newblock A gradient structure for reaction-diffusion systems and for
  energy-drift-diffusion systems.
\newblock {\em Nonlinearity}, 24:1329--1346, 2011.

\bibitem{Miel133TMER}
A.~Mielke.
\newblock Thermomechanical modeling of energy-reaction-diffusion systems,
  including bulk-interface interactions.
\newblock {\em Discr. Cont. Dynam. Systems Ser.~A}, 6:479--499, 2013.

\bibitem{MieRou15RIST}
A.~Mielke and T.~Roub{\'\i}{\v{c}}ek.
\newblock {\em Rate-Independent Systems -- Theory and Application}.
\newblock Springer, New York, 2015.

\bibitem{Mind65SGSS}
R.~Mindlin.
\newblock Second gradient of strain and surface-tension in linear elasticity.
\newblock {\em Intl. J. Solids Structures}, 1:417--438, 1965.

\bibitem{Pawl06TCCH}
I.~Paw{\l}ow.
\newblock Thermodynamically consistent {C}ahn-{H}illiard and {A}llen-{C}ahn
  models in elastic solids.
\newblock {\em Disc. Cont. Dynam. Syst.}, 15:1169--1191, 2006.

\bibitem{PawZaj08WSCH}
I.~Paw{\l}ow and W.~Zajaczkowski.
\newblock Weak solutions to 3-{D} {C}ahn-{H}illiard system in elastic solids.
\newblock {\em Topological Meth. Nonlin. Anal.}, 32:347--377, 2008.

\bibitem{Roub17VMSS}
T.~Roub{\'{\i}}{\v{c}}ek.
\newblock Variational methods for steady-state {D}arcy/{F}ick flow in swollen
  and poroelastic solids.
\newblock {\em Zeit. angew. Math. Mech.}, 97:990--1002, 2017.

\bibitem{RouSte18TEPR}
T.~Roub{\'{\i}}{\v{c}}ek and U.~Stefanelli.
\newblock Thermodynamics of elastoplastic porous rocks at large strains towards
  earthquake modeling.
\newblock {\em SIAM J. Appl. Math.}, 78:2597--2625, 2018.

\bibitem{RouTom13TMPE}
T.~Roub{\'{\i}}{\v{c}}ek and G.~Tomassetti.
\newblock A thermodynamically consistent model of magneto-elastic materials
  under diffusion at large strains and its analysis.
\newblock {\em Zeit. angew. Math. Phys.}, 64:1--28, 2013.

\bibitem{RouTom14THSM}
T.~Roub{\'{\i}}{\v{c}}ek and G.~Tomassetti.
\newblock Thermomechanics of hydrogen storage in metallic hydrides: modeling
  and analysis.
\newblock {\em Discrete Cont. Dynam. Systems B}, 19:2313--2333, 2014.

\bibitem{Toup62EMCS}
R.~Toupin.
\newblock Elastic materials with couple stresses.
\newblock {\em Arch. Ration. Mech. Anal.}, 11:385--414, 1962.

\bibitem{WPMB14GELD}
T.~Waffenschmidt, C.~Polindara, A.~Menzel, and S.~Blanco.
\newblock A gradient-enhanced large-deformation continuum damage model for
  fibre-reinforced materials.
\newblock {\em Comput. Methods Appl. Mech. Engrg.}, 268:801--842, 2014.

\bibitem{YaBa07TCMS}
V.~Yashin and A.~Balazs.
\newblock Theoretical and computational modeling of self-oscillating polymer
  gels.
\newblock {\em J. Chem. Phys.}, 26:Art.no.124707, 2007.

\bibitem{YSYB12CDBH}
V.~Yashin, S.~Suzuki, R.~Yoshida, and A.~Balazs.
\newblock Controlling the dynamic behavior of heterogeneous self-oscillating
  gels.
\newblock {\em J. Mater. Chem.}, 22:Art.no.13625, 2012.

\end{thebibliography}
%








\COMMENT{TO CHECK:
P. Steinmann, Formulation and computation of geometrically non-linear gradient damage, International Journal for Numerical Methods in
Engineering 46 (1999) 757--779}

\COMMENT{TO CHECK:    Richard Ostwald Ellen Kuhl Andreas Menzel
  On the implementation of finite deformation gradient-enhanced damage models
  Comput Mech (2019). https://doi.org/10.1007/s00466-019-01684-5}

\end{sloppypar}
\end{document}